\input amstex
\documentstyle{amsppt}
\document
\magnification=1200
\NoBlackBoxes
\nologo
\hoffset=0.7in
\voffset=1in
\def\wR{\widehat{R}}

\def\Q{\bold{Q}}
\def\Z{\bold{Z}}
\vsize16cm


\bigskip

 \centerline{\bf CYCLOTOMY}

\medskip

\centerline{\bf  AND ANALYTIC  GEOMETRY OVER $F_1$}

\medskip

\centerline{\bf Yuri I. Manin}

\medskip

\centerline{\it Max--Planck--Institut f\"ur Mathematik, Bonn, Germany,}

\centerline{\it and Northwestern University, Evanston, USA}

\bigskip

{\bf Abstract.}  Geometry over  non--existent  ``field with one element''  $F_1$ conceived 
by Jacques Tits [Ti]  half a century ago recently found an incarnation,
in several related but different guises.
In this paper I analyze the crucial role of roots of unity in this geometry
and propose a version of the notion of ``analytic functions'' over $F_1$.
The paper combines a focused survey of various approaches with some new constructions.

\bigskip

\hfill{\it To Alain Connes, for his sixtieth anniversary}

\medskip

\centerline{\bf 0. Introduction: many faces of cyclotomy}

\medskip

{\bf 0.1. Roots of unity and  field with one element.} 
The basics of algebraic geometry over
an elusive ``field with one element $F_1$"  were laid down recently in [So], [De1], [De2], [TV],
fifty years after a seminal remark by J.~Tits [Ti]. There are many motivations to
look for $F_1$; a hope to imitate Weil's proof for Riemann's zeta is one of them,
cf. [CCMa3], [Ku], [Ma1].

\smallskip

An important role in the formalization of $F_1$--geometry was played by the suggestion made in [KS] that one should simultaneously
consider all the ``finite extensions'' $F_{1^n}.$ This resulted in the approach of [So], where
a geometric object, say a scheme, $V$ over $F_1$, acquired flesh after a base extension to $\bold{Z}$,
and the $F_1$-- geometry of $V$ was reflected in (and in fact, formally defined in terms of) 
the geometry of  ``cyclotomic''  points of  an appropriate ordinary scheme $V_\bold{Z}$.  
In [De1] and [TV],  schemes over $F_1$ are defined  in categorical terms
independently of cyclotomy, but the latter reappears soon: see the Definition 1.7.1 below and the following discussion.

\smallskip 

All these ideas are interrelated but lead to somewhat different versions of basic definitions,
and develop  the initial  intuition in different directions,  so that their divergence can be
fruitfully exploited.  With this goal in mind, I have chosen the topics to be discussed
in sec. 1, where four approaches to the definition of $F_1$--geometry are sketched and compared.

\smallskip

Of course, roots of unity appear naturally  in many different geometric contexts, not motivated
by  geometry over $F_1$: some of  these contexts are  reviewed below in the subsections 0.2--0.6
of this Introduction. 
I have compiled a sample of them with an explicit goal: to guess how the insights gained within these contexts could
help develop $F_1$--geometry. 

\smallskip

Seemingly, a similar desire moved the authors of [CCMa2] to
put the theory of Bost--Connes in the framework of $F_1$--geometry.

\smallskip

I show in sec. 2 and sec. 3
that the results of [Ha1], the preparatory part to K.~Habiro's work
on Witten--Reshetikhin--Turaev invariants of homology spheres [Ha2],
and discoveries about these invariants made in [Law] and [LawZ], can be naturally viewed
as a contribution to the rudiments of analytic geometry over $F_1$.

\smallskip

Finally, in sec. 4 I discuss Witt vectors and a series of $F_1$--models of moduli spaces.

\smallskip

{\it Acknowledgements.} V.~Golyshev's  note [Go] prompted me 
to think about cyclotomy in the $F_1$ context. 
K.~Habiro read a preliminary version of this paper and suggested several
complements and simplifications. H.~Lenstra kindly refered me to [van D]
and other useful sources on profinite numbers.  A.~Connes and C.~Consani
sent me a copy of their new paper [CC] which was being written during the same
weeks as the first version this article. 

\smallskip

When this first version  appeared in arXiv
on Sept. 09, 2008, I received several messages 
commenting upon and developing the framework
involving roots of unity in $F_1$--geometry. 

\smallskip

Matilde Marcolli  used the multivariable Habiro ring in [Marc] in order to
generalize the Bost--Connes system.

\smallskip

James Borger drew my attention to the fact that my treatment 
of the cyclotomic coordinates on Witt schemes perfectly 
matches his remarkable  basic
idea that  ``a lambda--ring structure (in the sense of Grothendieck--Riemann--Roch) on a ring $R$ should be thought of as descent data
for $R$ from $\Z$ to $F_1$'' (message of Sept. 11, 2008).
Borger's approach promises to be a significant breakthrough in our understanding
of $F_1$--geometry, and I have added a brief discussion of it 
in this new version.

\smallskip

Finally, a totally anonymous referee provided a list of useful remarks and suggestions.

\smallskip
I  deeply appreciate 
their interest and help.

\medskip

{\bf 0.2. Roots of unity and Morse--Smale diffeomorphisms.} This aspect of cyclotomy
is described by D.~Grayson in [Gr]. 

\smallskip

Let $M$ be a compact smooth manifold, $f$ a diffeomorphism of $M$.
It is called {\it Morse--Smale}, if it is structurally stable,  and only a finite
number of points $x$ are non--wandering. (A point $x$ is called
non--wandering, if  for any neighborhood $U$ of $x$,
we have $U\cap f^n(U)\ne \emptyset$ for some $n>0$).

\smallskip

Assume that all eigenvalues of the action of $f$ on integral cohomology of $M$ are roots of unity
and put the question: {\it when $f$ is isotopic to a Morse--Smale map?}

\smallskip

There is an obstruction to this, lying in the group $SK_1(\Cal{R})$, where $\Cal{R}$ is the ring obtained by
localizing $\bold{Z}[q]$ with respect to $\Phi_0(q):=q$ and all cyclotomic polynomials
$$
\Phi_n(q):=\prod_{\eta} (q-\eta )
$$
where $\eta$ runs over all primitive roots of unity of degree $n\ge 1$.

\smallskip

This ring turns out to be a principal ideal domain. The reason for this is that
each closed point (a prime ideal of depth two) of the ``arithmetical plane''  $Spec\,\bold{Z}[q] $
is situated on an arithmetic curve $\Phi_n(q)=0, n\ge 0,$ because all finite fields
consist of roots of unity and zero.

\smallskip
Localization cuts all these curves off,
 and all  closed points go with them. The remaining prime ideals are of height one, and they are
 principal.
 
 \smallskip
 
 The same effect can be achieved by localizing wrt all primes $p\in \bold{Z}$, thus getting 
 the principal ideal domain $\bold{Q}[q]$.
 This localization  cuts away the closed fibers of the projection $Spec\,\bold{Z}[q] \to Spec\,\bold{Z}$,
 and all the closed points with them.

\smallskip

This suggests that the union of cyclotomic arithmetic curves
$\Phi_n(q)=0$  can be imagined as the union of closed fibers
of the projection $Spec\,\bold{Z}[q] \to Spec\,F_1[q]$, and the arithmetic plane itself
as the product of two coordinate axes, arithmetic one $Spec\,\bold{Z}$
and geometric one, $Spec\, F_1[q]$, over the ``absolute point'' $Spec\, {F}_1$.

\smallskip

In sec. 1 below, I review several versions of algebraic $F_1$--geometry where this intuition can
be made precise.

\smallskip

{\bf 0.2.1.} {\it Question.} Is there a context in which diffeomorphisms $f$, acting on
integral cohomology of $M$ with eigenvalues roots of unity, could be interpreted
as  ``Frobenius maps in caracteristic 1'', and their fixed (or non--wandering) points
in a Morse--Smale situation as
$F_{1^n}$--points of an appropriate variety?

\medskip

{\bf 0.3. Roots of unity and the Witten--Reshetikhin--Turaev invariants.}
An apparently totally different line of thought led to the consideration of completions
of $\bold{Z}[q]$ with respect to various  linear topologies generated by
the cyclotomic polynomials $\Phi_n(q)$.  Namely, it turned out that
the invariants of 3--dimensional homology spheres,
introduced first by E.~Witten by means of path integrals, and mathematically
constructed by Reshetikhin and Turaev, can be unified
into objects lying in completions of the kind described above.

\medskip

{\bf 0.3.1.} {\it  Question.} Can these completions be interpreted in a framework
of $F_1$--geometry?

\smallskip

We try to answer this question affirmatively in sec. 2  and 3 below.

\smallskip
(Similar completions along
 the arithmetical axis produce for example direct products of $p$--adic integers
$\prod \bold{Z}_{p_i}$ and the ring $\varprojlim_N\bold{Z}/(p_1\dots p_N)$,
in which $\bold{Z}$ can be embedded.)

\smallskip

   We suggest two interpretations, one in the framework of Soul\'e's axiomatics,
and another more in the spirit of To\"en--Vaqui\'e and Deitmar's definitions.
Here are brief explanations.

\smallskip

Soul\'e's   definition of an $F_1$--scheme $X$ involves, besides $X_{\bold{Z}}$,
a $\bold{C}$--algebra $\Cal{A}_X$, and each cyclotomic point
of $X_{\bold{Z}}$ coming from $X$ must assign ``values''  to the elements
of $\Cal{A}_X$. His choice of $\Cal{A}_X$ for the multiplicative group
$G_{m,F_1}$ is that of continuous functions on the unit circle in $\bold{C}$
(cf. [So], 5.2.2).  For the affine line he uses
holomorhic functions in the open unit circle continuous
on the boundary.

\smallskip

We suggest to consider respectively  {\it the ring of
Habiro's analytic functions} and {\it the ring of Habiro's functions
admitting an analytic continuation in the open  unit disc.}
The first one consists of formal series
$$
f(q) = a_0 +\sum_{n=1}^{\infty} a_n(1-q)\dots (1-q^n).
\eqno(0.1)
$$
where $a_n$ are polynomials in $q$ of degree $\le n-1.$
At any root of unity, only a finite number of terms do not vanish,
so $f$ is a well  defined function on cyclotomic points.

\smallskip

The second option consists in considering holomorphic functions $\varphi (q)$ in the unit circle,
such that for any root of unity $\zeta$, a radial limit $\lim_{r\to 1_{-}} r\zeta$ exists,
and the family of such limits can be given by the series (0.1).

\smallskip

Versions of this choice might involve functions holomorphic inside
variable rings with outer boundary $|q|=1$ admitting radial limits
at roots of unity, or even pairs  of functions $\varphi_{-}$, resp. $ \varphi_{+}$, holomorphic inside narrow rings
with outer, resp. inner, boundary unit circle,  which restrict to a $C^{\infty}$--function on all
small radial intervals  $(1-\varepsilon_{\zeta}, 1+\varepsilon_{\zeta})\cdot \zeta$
containg roots of unity $\zeta$. In particular, they must satisfy

$$
\lim_{r\to 1_{-}} \varphi_{-}\,( r\zeta ) = \lim_{r\to 1_{+}} \varphi_{+} ( r\zeta )
$$
The limit values should admit the representation (0.1).

\smallskip

The fact that there exist highly nontrivial and interesting examples of such functions,
was discovered in the theory of Witten--Reshetikhin--Turaev invariants: cf. [Law], [LawZ]. Don Zagier
   says that $\varphi_{\pm}$ ``leak through'' roots of unity.

\smallskip

On the other hand, if $\Cal{A}_X$ is not a part of the definition
of a $F_1$--scheme, as in the versions of [TV] and [De1],
one can still imagine that a ring of the type discussed above
would form a part of the structure of analytic $F_1$--varieties
when $F_1$--geometry becomes mature enough to include
analytic geometry.

\smallskip

[CCMa2] also suggests that  time is ripe for such generalizations.

\medskip

{\bf 0.4. Roots of unity and the Bost--Connes system.}
In the paper [BoCo] roots of unity appear in the following setting. Consider the Hecke algebra $\Cal{H}$
with involution  over 
$\bold{Q}$
given by the following presentation. The generators
are denoted $\mu_n,\,n\in \bold{Z}_+$, and $e(\gamma )$, $\gamma\in
\bold{Q}/\bold{Z}.$ The relations are
$$
\mu_n^*\mu_n=1,\quad \mu_{mn}=\mu_m\mu_n,\quad
\mu_m^*\mu_n=\mu_n\mu_m^*\ \roman{for}\ (m,n)=1;
$$
$$
e(\gamma)^*=e(-\gamma ),\quad e(\gamma_1+\gamma_2) =
e(\gamma_1)\,e(\gamma_2);
$$
$$
e(\gamma )\mu_n= \mu_n e(n\gamma ),\quad 
\mu_n e(\gamma )\mu_n^*=\frac{1}{n}
\sum_{n\delta =\gamma} e(\delta ).
$$
The id\`ele class group $\widehat{\bold{Z}}^*$ of $\bold{Q}$ acts upon $\Cal{H}$
in a very explicit and simple way: on $e(\gamma )$'s
the action is induced by the multiplication $\widehat{\bold{Z}}^*\times \bold{Q}/\bold{Z}
\to \bold{Q}/\bold{Z}$, whereas on $\mu_n$'s it is identical.

\smallskip

The algebra $\Cal{H}$ admits an involutive representation $\rho$ in $l^2
(\bold{Z}_+)$: denoting by $\{\epsilon_k\}$ the standard basis
of this space, we have
$$
\rho (\mu_n)\,\epsilon_k = \epsilon_{nk},\quad
\rho (e(\gamma ))\epsilon_k= e^{2\pi ik\gamma}\epsilon_k.
$$
From this, one can produce the whole
$\roman{Gal}\,(\bold{Q}^{ab}/\bold{Q})$--orbit $\{\rho_g\}$ of 
such representations,
applying $g\in \roman{Gal}\,(\bold{Q}^{ab}/\bold{Q})$
to all roots of unity occuring at the right hand sides of
the expressions for  $\rho (e(\gamma ))\epsilon_k$.
All these representations can be canonically extended
to the $C^*$--algebra completion $C$ of $\Cal{H}$
constructed from the regular representation of $\Cal{H}$.
Let us denote them by the same symbol $\rho_g$.

\smallskip

To formulate the main theorem of [BoCo], we need some more explanations.
The algebra $C$ admits a  canonical action of $\bold{R}$,
which can be interpreted as {\it time evolution} 
represented on the algebra of observables. This is
a general (and deep) fact in the theory of $C^*$--algebras,
but for $C$ the action of $\bold{R}$ can be quite explicitly
described on the generators. Let us denote by $\sigma_t$
the action of $t\in\bold{R}.$ {\it A KMS${}_{\beta}$ state
at inverse temperature} $\beta$ on $(C,\sigma_t)$
is defined as a state $\varphi$ on $C$ such that for any
$x,y\in C$ there exists a bounded holomorphic function
$F_{x,y}(z)$ defined in the strip $0\le \roman{Im}\,z\le
\beta$ and continuous on the boundary, satisfying
$$
\varphi (x\sigma_t(y)) = F_{x,y}(t),\ 
\varphi (\sigma_t(y)x) = F_{x,y}(t+i\beta ).
$$
Now denote by $H$ the positive operator on
$l^2(\bold{Z}_+)$: $H\epsilon_k = (\roman{log}\,k)\,\epsilon_k.$
Then for any $\beta >1$, $g\in \roman{Gal}\,(\bold{Q}^{ab}/\bold{Q})$
one can define a KMS${}_{\beta}$ state $\varphi_{\beta,g}$
on $(C,\sigma_t)$ by the following formula:
$$
\varphi_{\beta,g}(x):=\zeta (\beta )^{-1}\,\roman{Trace}\,
(\rho_g(x)\,e^{-\beta H}),\quad x\in C
$$
where $\zeta$ is the Riemann zeta--function.
The map $g\mapsto \varphi_{\beta,g}(x)$ is a homeomorphism
of $\roman{Gal}\,(\bold{Q}^{ab}/\bold{Q})$ with the space of extreme
points of the Choquet simplex of all KMS${}_{\beta}$ states.

\smallskip

To the contrary, for $\beta <1$ there is a unique KMS${}_{\beta}$ state.
This is a remarkable ``arithmetical symmetry breaking'' phenomenon.

\smallskip

The description of the Hecke algebra above involves denominators
in the last relation. In [CCMa2], the authors construct $\bold{Z}$--models
of finite layers of this object and natural morphisms between
them, and show that the resulting system is a lift to $\bold{Z}$
of an $F_1$--tower.

\smallskip

This picture is generalized to the multivariable case in [Marc].

\medskip

{\bf 0.5. Witt vectors.} It is desirable to consider
the arithmetical axis $Spec\,\bold{Z}$ 
as an $F_1$--space as well, but in the current framework
it is certainly not a scheme of finite type. In fact,
its base extension to $\bold{Z}$ is elusive, being precisely
what we would like to see as the spectrum of $\bold{Z}\otimes\bold{Z}$.

\smallskip

Nevertheless, in a certain sense  primes can be considered
as cyclotomic points of $Spec\,\bold{Z}$, at which the ``cyclotomic coordinates'',
all integers, take values that are roots of unity or zero.

\smallskip

In fact, roots of unity of degree $q=p^{n}-1$ (and zero), considered
together with their embedding into  a fixed unramified extension  $\bold{Z}_p^{nr}$ 
of $\bold{Z}_p$
rather than $\bold{C}$, appear as natural coefficients of $p$--adic expansions
discovered by Teichm\"uller and Witt. Namely, each residue class  in $\bold{Z}_p^{nr}/(p)$
has a unique (Teichm\"uller) representative $\zeta$ which is either a root of unity or $0$ in $\bold{Z}_p^{nr}$,
so that an element of such an extension can be written as a well defined series
$\sum_{i=0}^{\infty} \zeta_ip^i$. Moreover, coefficients of a sum or a product
of two such series are given by Witt's universal polynomials in the coefficients
of the summands/factors in the following sense: one must reduce Witt's coefficients 
modulo $p$, apply these polynomials (which are defined over $\bold{Z}$),
and lift results back to roots of unity.

\smallskip

This can be generalized to the so called
``big Witt ring'' and interpreted in the following way. On affine spaces
$A_{\bold{Z}}^{k} =Spec \,\bold{Z}\,[u_1, u_2, \dots , u_k]$ there exists a natural system
of ``cyclotomic coordinates'' (in the $p$--adic context sometimes called ``ghost coordinates'').
In terms of these coordinates, one can define an $F_1$--gadget  \`a la Soul\'e,
requiring that in the subfunctor of points, these coordinates took cyclotomic values
(including  zero). However, Witt's addition/multiplication becomes well defined
only after extension to $\bold{Z}$, unless the notion of morphism 
over $F_1$ is drastically extended. 

\smallskip

To me, this looks like  a strong argument for considering
options for such an extension.

\smallskip

We supply some more details in sec. 4. 

\medskip

{\bf 0.6. Roots of unity and the analogy between Hilbert polynomials and zeta functions.}
There were several suggestions  that Hilbert polynomials $H(n)$, say, of graded commutative
rings behave like toy zeta functions. 

\smallskip

Rather precise recent observations by V.~Golyshev in [Go]  can be summarized as
follows.  

\smallskip

a) The comparison to zetas becomes most striking if one restricts oneself to the following
Hilbert polynomials of projective smooth manifolds $X$:

\smallskip

(i) If $X$ is of general type or Fano: consider $H_{-K_{X}}(n) := \chi (-nK_X)$.

\smallskip

(ii) If $X$ is a Calabi--Yau manifold embedded as an anticanonical section in a Fano manifold:
consider the Euler characteristic of the powers of the induced anticanonical sheaf. 

\smallskip

b) With this normalization, the Serre duality leads to a functional equation  for the Hilbert
polynomial of the $s\mapsto -1-s$ type.

\smallskip

c) In many cases, the well known inequalities for Chern numbers of $X$ 
imply that all roots of $H(s)$ lie in the critical strip $-1 < Re\,s < 0,$
and sometimes even more precise statements. 
For example, Yau's inequality $c_1^3\ge 8/3 c_1c_2$ for Fano's threefolds
shows this fact for them.

\medskip

{\bf 0.6.1} {\it Question.}  Is there a systematic relationship between Hilbert polynomials
and zeta functions of schemes  (or more general spaces) over $F_1$?

\smallskip

The existence of such zeta functions and their structure
in certain cases was heuristically suggested in [Ma1]
(cf. also [Ku]). They make precise sense  for some specimens in Soul\'e's
category, and are indeed polynomials; see also [CC] for essential
complements. An obvious  attack on question 0.6.1 might 
start with comparing the counting of $F_{1^n}$--points with counting of
monomials in cyclotomic coordinates.

\smallskip

Roots of unity appear in this context  via the following beautiful observation due to
F.~Rodriguez--Villegas, [RV].

\smallskip

Consider first an arbitrary polynomial $H(q)\in \bold{Z}[q]$. Define another polynomial
$P(t)$ such that 
$$
\sum_{n=0}^{\infty} H(n)t^n= \frac{P(t)}{(1-t)^d}
$$
where  $P(1)\ne 0$. Let $d:=\roman{deg}\,H+1,\ e:= \roman{deg}\, P.$

\smallskip

Rodriguez--Villegas proves that {\it if all roots of $P$ lie on the unit circle,
then $H(q)$ has simple roots at $q=-1,\dots, e+1-d$ and possible additional roots 
at the middle of this critical strip $Re \, q=\dfrac{e-d}{2}.$}

\smallskip

This result can be applied to the case when $H(n)=\roman{dim}\, A_n$ is the Hilbert
polynomial of a graded algebra 
$\oplus_{n=0}^{\infty} A_n$ generated by $A_1$, of Krull dimension $d$,
complete intersection of polydegree $(n_1,\dots ,n_s)$.

\smallskip

It turns out that $e=n_1+\dots +n_s$ and $P(t)=\prod_{j=1}^s (1+t+\dots +t^{n_j})$
so that all roots of $P$ are in fact roots of unity.

\smallskip

The critical strip in [Go] has width 1, because Golyshev differently normalizes
the grading via $-K_X$. However, the Rodriguez--Villegas grading agrees with the
motivic philosophy involving weights and Tate's motives over $F_1$, see [Ma1].

\medskip

{\bf 0.7. Summary.} We are guided by the following heuristics. Each time
that roots of unity appear in a certain context, we try to interpret
the functions whose values are these roots of unity as 
cyclotomic coordinates on a relevant $F_1$--scheme, in the sense 
of the Definition 1.7.1 below, or a version
thereof.

\smallskip

An appropriate version of (big) Witt vectors must furnish the basic 
$F_1$--analytic (or formal) approximation to the arithmetic line $Spec\,\bold{Z}$.

\bigskip
\centerline{\bf 1. Geometry over $F_1$: generalities}

\medskip

This section sketches and compares  four approaches to the definition of $F_1$--geometry.
Preparing a colloquium talk in Paris, I have succumbed to the
temptation to associate them with some dominant trends in the history of art. 

\medskip
{\bf 1.1. Affine schemes over $F_1$ according to To\"en  and Vaqui\'e}
{\it (Abstract Expressionism).} Affine schemes over $F_1$ arise in the most straightforward (and allowing vast generalizations) manner  in the framework 
of  [TV], according to which
algebraic geometry over $F_1$ is a special case of algebraic geometry relative to
 a monoidal symmetric category $(C,\otimes ,\bold{1})$, which is assumed to be
complete, cocomplete, and to admit internal $Hom$'s.

\smallskip

Such a category $C$ gives rise to the category of commutative, associative and unitary monoids
$Comm\, (C)$ which serves as a substitute for the category of ordinary commutative rings.
Each object $A$ of $Comm\,(C)$ determines the category of $A$--modules $A$--$Mod$
consisting of pairs  $(M,\mu)$ where $M$ is an object of $C$ together with action $\mu:\,A\otimes M\to M$ 
and satisfying the usual formalism.

\smallskip

The opposite category $Aff_C:=Comm\,(C)^{opp}$ is called the category of
$C$--affine schemes, and the tautological functor $Comm\,(C)\to Aff_C$ is called $Spec.$

\smallskip

Florian Marty in [Mart2] defines and studies the notion of smoothness
in the To\"en--Vaqui\'e geometry. This requires passing to the homotopical
algebra in appropriate simplicial categories.

\smallskip 

 According to [TV], we obtain $F_1$--geometry as the geometry
relative to the monoidal category of sets and direct products
$(Ens ,\times ,*)$.  

\smallskip

Commutative rings relative to   $(Ens ,\times ,*)$ are just 
the ordinary commutative (associative, unital) monoids written multiplicatively:
this explains the popular motto that to do $F_1$--algebra one must forget
the additive structure: cf. [Har]. This structure is restored when one
applies the functor ``base change'' $\otimes_{F_1}\bold{Z}$: a monoid $M$
turns into the commutative associative unital ring $\bold{Z}[M]$. 
The opposite to monoids category will be denoted $Aff_{F_1}$.

\smallskip

More generally, any commutative ring $R$ determines the base extension functor
$$
\otimes_{F_1}R:\,\, Aff_{F_1}\to Aff_R,\
M\mapsto R[M],
$$
from affine schemes over $F_1$ to affine schemes over $Spec\, R$.

\smallskip

Elements of monoids $M$ will be called {\it cyclotomic coordinates}
on the respective affine scheme. The same term will refer to
their images in $R[M]$. On more general schemes, we may speak about
local cyclotomic coordinates.

\smallskip

{\bf 1.2. Deitmar's affine schemes} {\it (Minimalism).} A.~Deitmar in [De1] adopts the same definition of the category $Aff_{F_1}$.
Moreover, he associates to a monoid $M$ a topological space
which we will denote $spec\,M$ (to distinguish it from the spectrum of prime ideals of a ring $Spec\,(*)$), and which is endowed 
with a structure sheaf. Points of this space are prime ideals $P\subset M$: such submonoids that
$xy\in P$ implies $x\in P$ or $y\in P$.
Basic open sets and the structure (pre)sheaf of monoids are determined via localization, 
just as in the classical case of commutative rings.
Moreover, Deitmar characterizes morphisms in $Aff_{F_1}$ in terms of appropriate
morphisms of topological spaces $spec$ with structure sheaves.

\medskip

{\bf 1.3. Examples.} {\it (i) Affine $F_1$--schemes associated to abelian groups.} Let $M$ be an
abelian group considered as a monoid in $Ens$. We have
$$
spec\,M\otimes_{F_1}\bold{Z} = Spec\,\bold{Z}[M] .
$$
In particular, [TV] define $F_{1^n}$ as the monoid (group) $\bold{Z}/n\bold{Z}$,
and its spectrum after lifting to $\bold{Z}$ becomes 
$$
spec\, F_{1^n}\otimes_{F_1}\bold{Z}:= Spec\, \bold{Z}[q]/(q^n-1)= Spec\, \bold{Z}[q]/\{n\}_q.
\eqno(1.1)
$$
In [So], the study of $\bold{Z}[q]/\{n\}_q$--points of an a priori given ordinary scheme $X$
gives clues to finding its $F_1$--forms identified with certain subfunctors
of $F_{1^n}$--points.
\smallskip

In our paper, formula (1.1) motivates the introduction of analytic functions
on (certain) $F_1$--schemes via Habiro's formalism:   morally, they are functions that are defined 
{\it at all $F_{1^n}$--points, but nowhere else.}  (In fact, the latter  stricture should not be
taken too literally: some functions  have very interesting $p$--adic,
and sometimes complex, arguments and values as well).

\medskip

{\it (ii). Affine scheme $G_{m,F_1}$.} Over $F_1$, it is represented by the spectrum of 
the infinite cyclic group $\bold{Z}$. Lifted to $\bold{Z}$, it becomes
the ordinary $G_m=Spec\,\bold{Z}[q,q^{-1}].$ 

\smallskip


\medskip

{\bf 1.4. Affine spaces.} Affine line $A^1_{F_1}$ is the  spectrum of the infinite
cyclic monoid $\bold{N}$. Its lift to $\bold{Z}$ is $A^1_{\bold{Z}}:=Spec\,\bold{Z}[q].$
Similarly, $A^k_{F_1}$ is the spectrum of $\bold{N}^k$, $k\ge 1.$

\smallskip
The space $spec\,\bold{N}$ consists of one closed point $(q)$ and one generic point.

\smallskip

One can also consider $\bold{N}^{\times}$ that is, the free monoid freely generated
by all primes.  Its lift to $\bold{Z}$ is the ring of polynomials
in infinitely many variables indexed by primes.

\medskip

{\bf 1.5. Affine scheme $GL(n)_{F_1}$}. According to [TV], Proposition 4.1,  
the natural sheaf (in the Grothendieck  flat topology, see [TV] and 1.6 below)
of automorphisms of a free module of rank $n$ is represented,
after lifting to $\bold{Z}$, by the semidirect product of $G_m^n$ and $\bold{S}_n$
(the symmetric group). 

\smallskip

In the more down--to--earth language of [De1], sec. 5,
this is expressed as follows. Let $A$ be a commutative monoid.
Define the (set--theoretic) group of ``$A$--valued points of $GL(n)_{F_1}$''
as
$$
GL(n)_{F_1}(A):= \roman{Aut}_A (A^n).
$$
This  can be identified with the group of $(n,n)$ matrices with entries
in $A\subset \bold{Z}[A]$, having exactly one non--zero element in each
row and each column.  This is precisely the description
of [TV] quoted above.

\smallskip

The reader should be warned that, unlike to what happened with $A^k$ and $G_m$,   after lifting $GL(n)_{F_1}$ to $\bold{Z}$
we {\it do not} get the usual $GL(n)_{\bold{Z}}$.  This caused a difficulty in the framework
of [So], where it was not obvious how to choose ``cyclotomic points of $GL(n)_{\bold{Z}}.$''
In fact, according to [TV], Proposition 4.1, $GL(n)_{\bold{Z}}$ for $n>1$ is not
a lift of an $F_1$--scheme in their sense.

\medskip

{\bf 1.6. General schemes over $F_1$.} Glueing general schemes from affine ones 
is defined differently in [TV] and [De1] respectively.

\smallskip
For Deitmar, an $F_1$--scheme is a topological space with a sheaf of monoids
that is everywhere locally affine, that is, locally isomorphic to some $spec\, M$.

\smallskip

To\"en and Vaqui\'e endow the category $Aff_C$ with a natural Grothendieck's topology, which is called {\it the flat topology.}
Using it, one can defined general {\it schemes} relative to $C$, as functors that can be obtained
from disjoint unions of affine schemes $X$ by taking the quotient with respect to
an equivalence relation $R\subset X\times X$ such that projections $R\to X$
are local Zariski isomorphisms. Such schemes form a category denoted $Sch\,(C)$.

\smallskip

Florian Marty in [Mart1] presents a thorough study of Zariski topology on the category
of commutative monoids in $C$ and applies it to the comparison
of Deitmar's schemes with To\'en--Vaqui\'e's ones.

\medskip

{\bf 1.7. Schemes over $F_1$ \`a la Soul\'e} {\it (Critical Realism).} The idea of Soul\'e's definitions in [So] can be succintly 
formulated as the project  of direct reconstruction 
of $F_1$--schemas $X$ {\it of finite type} from certain schemes $X_{\bold{Z}}$ over $\bold{Z}$
endowed with some kind of {\it descent data} from $\Z$ to $F_1$.

\smallskip

However, more than only descent data to $F_1$ is required: Soul\'e's spaces   
come with an additional
data $\Cal{A}_X$ which is a $\bold{C}$--algebra, morally an algebra of functions
on the ``$\infty$--adic  completion'' of $X$.

\smallskip

This latter structure embeds $F_1$--geometry into a wider context, potentially containing also rich structures of
Arakelov, or $\infty$--adic geometry. Some hints that this should be necessary and possible
can be glimpsed in the remark made in [Ma1], 1.7. Namely, in [Ma1] it was suggested
that the zeta function of $\bold{P}_{F_1}^k$ must be $(2\pi )^{-(k+1)}s(s-1)\dots (s-k).$

\smallskip
Combining this with  Deninger's representation of the basic Euler $\Gamma$--factor
at arithmetical infinity as a regularized product
$$
\Gamma_{\bold{C}}(s)^{-1}:= \frac{(2\pi)^s}{\Gamma (s)} =\prod_{n\ge 0} \frac{s+n}{2\pi}
$$
we see that this gamma--factor should be understood as the zeta--function of the (motivic dual of)
an infinite dimensional projective space over $F_1$. 

\smallskip

However, the existing framework is too narrow to make sense of this statement:
although the zeta of $\bold{P}_{F_1}^k$ is now defined in [So] and agrees with
expectations of [Ma1] (up to a power of $2\pi$), the infinite--dimensional case and its connections
with $\infty$--adic geometry  still elude us. A promising approach extensively
elaborated in the thesis by N.~Durov (cf. [Du]) might pave the road to this unification.
The treatment of the Bost--Connes dynamical system in [CCMa2] 
provides another bridge between $F_1$--geometry and the archimedean world.

\smallskip

Returning to [So], we will now sketch his version of $F_1$--schemes of finite type.

\smallskip

The data defining such  a scheme $X$ consist of:

\smallskip

{ \it (i) A $\bold{Z}$--scheme of finite type $X_{\bold{Z}}$.

\smallskip

(ii) A subfunctor  $X(R)$ of the functor of points of $X_{\bold{Z}}$ from a category of rings
to the category of sets: 
$$
X_{\bold{Z}}(R):=\, Hom (Spec\,R,X_{\bold{Z}}).
$$}
Here $R$  runs over rings that are direct summands of  $\otimes_i \bold{Z}[q]/(\{n_i\}_q)$,
and  each $X(R)$ is required to be a finite set. We will call elements
of $X(R)$ ``cyclotomic points''.

\smallskip

{\it (iii) A $\bold{C}$--algebra $\Cal{A}_X$, and an assignement of complex values to each 
element  $f\in \Cal{A}_X$ at each pair consisting of  point of $X(R)$ and a ring
homomorphism $R\to \bold{C}$.}

\medskip

We will not spell out here the compatibility requirements between these data,
which are pretty straightforward. 

\smallskip

Morphisms of schemes over $F_1$ are pairs, consisting of functor
morphisms of cyclotomic points and contravariant homomorphisms
of function algebras, compatible with the rest of the data.

\smallskip

It is natural to call an $F_1$--scheme $X$ affine, if $X_{\bold{Z}}$ is affine. 
But without further restrictions, one would get many schemes over $\bold{Z}$
into which the cyclotomic points could be embedded as a subfunctor.
The restriction that restores the uniqueness of $X_{\bold{Z}}$ once $X$ is known
declares that $X_{\bold{Z}}$ must be the initial object in the category
of such embeddings (see [So], sec. 4, Definition 3, for a precise statement).
A  similar universality requirement defines general $F_1$--schemes
(loc. cit., Definition 5).

\smallskip

We will now formally define the notion of  {\it cyclotomic coordinates} on Soul\'e's
$F_1$--schemes. Let $X$ be an affine $F_1$--scheme, $X_{\bold{Z}}=Spec\,A$.

\medskip

{\bf 1.7.1. Definition.} {\it A cyclotomic coordinate on the affine $F_1$--scheme $X$
in the sense of Soul\'e is any element $f\in A$ whose values at all cyclotomic points $X(R)$ are either 0,
or roots of unity.}

\smallskip

Clearly, cyclotomic coordinates in this sense form a commutative monoid
with unit. If the scheme $X$ is not affine, local cyclotomic coordinates
can be defined, forming a (pre)sheaf of commutative monoids.

\smallskip

Recall that in the framework of [TV] and [De1], where $\bold{Z}$--lifts of $F_1$--spaces are patched from
spectra of monoid ring $\bold{Z}\,[S]$, the elements of $S$ themselves
were called cyclotomic coordinates. However, since these versions of $F_1$--schemes
are not equivalent to Soul\'e's one, we should use this term being aware of its context. 

\smallskip

Notice also that

\smallskip

a) To reconstruct cyclotomic coordinates in the sense 1.7.1, it is sufficient
to know $X_{\bold{Z}}$ and the functor $R\to X(R)\subset X_{\bold{Z}}(R).$ This is a part of the structure
of {\it a gadget}, as Soul\'e's {\it truc} was translated in [CCMa2]).

\smallskip

b) The rings $R$ used by Soul\'e to probe schemes over $F_1$ are essentially group rings
of finite abelian groups.

\smallskip

Conceivably, one could replace finite abelian groups by finite commutative unital monoids, thus
narrowing the gap between [So] and [TV], [De].

\smallskip

c) Moreover, one could sketch rudiments of  {\it supergeometry} over $F_1$,
by requiring $Z_2$--grading of our monoids, a structure subgroup
$\{\pm 1\}$, and the anticommutation rule for odd elements.

\smallskip

The following example from [So]  serves as a good illustration of similarities and differences
between affine schemes in the sense of [So] and [TV] respectively, and of relationships of $F_1$
to Arakelov geometry.

\medskip

{\bf 1.7.2. Arakelov vector bundles over $Spec\,\bold{Z}$ as affine $F_1$--schemes.}
An Arakelov vector bundle $\bar{\Lambda}$ over $Spec\,\bold{Z}$ is defined as a pair consisting
of a free abelian group $\Lambda$ of finite rank and an hermitean norm
$||\cdot ||$ on $\Lambda_{\bold{C}}:=\Lambda\otimes \bold{C}$, ``integral structure at arithmetical infinity''.
The global sections of $\bar{\Lambda}$ 
over the ``compactification''  $Spec\,{\bold{Z}}\cup \infty$ are defined as $B\cap \Lambda$, where
$B:= \{x\in\Lambda_{\bold{C}}\,|\,||x||\le 1\}.$

\smallskip

In order to produce a Soul\'e's affine scheme $X(\bar{\Lambda})$ out of $\bar{\Lambda}$, make an additional choice
(of which the final product will not depend):
choose a finite subset $\Phi\subset B\cap\Lambda\setminus \{0\}$ such that if
$v\in B\cap\Lambda\setminus \{0\}$, then exactly one element of the pair $\{v,-v\}$ belongs to $\Phi$.
Let $\Lambda_0$ be the sublattice of $\Lambda$ generated by $\Phi$, $\Lambda_0^t$ the dual lattice.

\smallskip

Now we can define the structure data.

\smallskip

(i) $X(\bar{\Lambda})_{\bold{Z}} :=\bold{Z}[\Lambda_0^t]$.

\smallskip
(ii) The  points of $X(\bar{\Lambda})(R)$ are given by the following
prescription:
$$
X(\bar{\Lambda})(R) := \{x=\sum_{v\in\Phi}v\otimes\zeta_v\,|\, x\in\Lambda\otimes_{\bold{Z}}R,\,\zeta_v\in \mu (R)\cup\{0\}\,\}.
$$
Equivalently, coefficients at $v\in \Phi$ are cyclotomic coordinates.

\smallskip

(iii) $\Cal{A}_{X(\bar{\Lambda})}$ is defined as the algebra of functions holomorphic
and continuous on the boundary of the following domain:
$$
C:=\{x\in \Lambda_0\otimes \bold{C}\,|\,||x|| \le \roman{card}\,\Phi\}.
$$
Given a homomorphism $\sigma :\,R\to \bold{C}$, a cyclotomic $R$--valued point $x$ and
a function $f\in \bold{Z}[\Lambda_0^t]$, we get its value at  $(x,\sigma )$ in an obvious way.

\smallskip

Comparing this example to the definitions of [TV] and [De1], we see
that  the algebra $X(\bar{\Lambda})_{\bold{Z}} :=\bold{Z}[\Lambda_0^t]$ fits in their framework, 
but that other elements
of the structure significantly change morphisms and points.

\medskip

{\bf 1.7.3. A non--affine case: toric varieties.}  The treatment of this case  in  [TV], [De1] and [So]
leads to essentially  the same object (although Soul\'e produces his $\bold{C}$--algebra only in
the smooth case i.~e. for regular fans).

\smallskip

Let $\Delta$ be a fan. Each element $\sigma\in \Delta$ determines the dual cone $\sigma^*$.
Let $M_{\sigma}$ be the commutative monoid of integer points of $\sigma^*$,
$U_{\sigma}$ is its spectrum as $F_1$--scheme.
If $\tau$ is a face of $\sigma$, we get a morphism of monoids $M_{\sigma}\to M_{\tau}.$ 
The respective morphism of schemes $U_{\tau}\to U_{\sigma}$ is Zariski open.
Put  $X:=\coprod_{\sigma\in \Delta} U_{\sigma}.$ According to [TV], 4.2,  the quotient of 
 $X\times X$ modulo equivalence relation $R:=\coprod_{\sigma ,\tau\in \Delta} U_{\sigma\cap\tau}$
defines an $F_1$--scheme $X(\Delta)_{F_1}.$ Lifting it to $\bold{Z}$,
we get the classical toric scheme $X(\Delta ).$

\smallskip

In [De1], the same quotient is straightforwardly interpreted as a glueing of   monoid spectra. 
In [So] the picture is enhanced by an appropriate $\bold{C}$--algebra.

\medskip

{\bf 1.8. Tits's problem and Connes--Consani schemes.} Tits remarked that 
one  can substitute $q=1$ in the classical formulas for the number of
$F_q$ points of a projective space $\bold{P}^{n-1}$ (resp.  Grassmanian
$Gr\,(n,j)$) and get  formulas for the cardinality of
$\{1,\dots ,n\}$ (resp. of the set of subsets of cardinality $j$ in it). 
Thus we get a version of classical combinatorial projective geometry, in which each line
has two points, each plane has three points etc. Tits asked in [Ti] how to extend
this to Chevalley groups and respective homogeneous spaces:
it would be a version of geometry  of homogeneous spaces ``over a field of characteristic 1''
as he put it then.

\smallskip

This project was realized only in 2008, when A.~Connes and C.~Consani
adapted Soul\'e's definition to this problem in [CC].  Their main innovation
consists in considering the functor of cyclotomic points $X(R)$ 
as taking values in the category of {\it graded
sets}. Only components of degree zero are taken in account in various
point counting contexts. After clarifying this issue  they find out
that Chevalley schemes have $F_{1^2}$ as a natural field of definition,
rather than $F_1$.

\medskip

{\bf 1.9. Lambda--rings and Borger's project}  {\it (Futurism).} As I have already mentioned,
the key idea of James Borger consists in a totally new conception of  $\Z$--{\it to}--$F_1$ descent data:
namely, {\it a  restricted $\lambda$--ring structure} in the sense of Grothendieck.  

\smallskip

According to [BorS], one can think about such a structure on a ring without additive torsion
$R$ as a family $\psi_p:\,R\to R$ of commuting ring endomorphisms indexed by primes such that
$\psi_p(x)-x\in pR$ for all $x, p$.

\smallskip
 
More generally, as is sketched in [Bor2], we may consider the category of ``spaces'' $\roman{Sp}_{\Z}$,
defined as sheaves of sets on the category of affine schemes with \'etale topology.
It is endowed with the endofunctor $W^*$ of infinite big Witt vectors (cf. the definition
in 4.1 below). This endofunctor carries a canonical monad structure.
A $\Lambda$--structure on a space $X$ is defined as an action of $W^*$ on $X$.
$\Lambda$--spaces with $W^*$--equivariant morphisms 
form a category $\roman{Sp}_{\Z /\Lambda}$. The functor forgetting
the $\Lambda$--structure is called $v^*:\, \roman{Sp}_{\Z /\Lambda} \to \roman{Sp}_{\Z}$.
It admits a left adjoint $v_!$ and a right one $v_*$.  The first one must be thought of
as (geometric) forgetting the base, and the second one as
Weil's restriction of scalars functor.

\smallskip

Using  the general topos formalism, Borger looks at {\it algebraic geometry of $\Lambda$--rings}
as a lifted {\it algebraic geometry over $\bold{F}_1$}, represented by
the big etale topos over $\bold{F}_1$.
 
\smallskip

 In particular,  the ring $W(\bold{Z})$ of big Witt vectors with entries in 
$\bold{Z}$
 should be thought of as (a completed version of)
 $\bold{Z}\otimes_{\bold{F}_1}\bold{Z }$.
 
\smallskip 

Varieties of finite type over $\bold{F}_1$ (in this sense) are very {\it rigid, combinatorial objects.}  They are essentially quotients of toric varieties by toric equivalence relations. 
In particular, only Tate motives descend to $\bold{F}_1$.

 \smallskip

 Non--finite--type schemes over 
 $\bold{F}_1$ are  more interesting.  The big de Rham--Witt cohomology of $X$ 
 ``is'' the de Rham cohomology of $X$ "viewed as an $\bold{F}_1$--scheme". It should contain the full information of the motive of $X$ and is probably a concrete universal Weil cohomology theory.  
 
 \smallskip

The Weil restriction of scalars from $\bold{Z}$ to $\bold{F}_1$ is an arithmetically global version of Buium's $p$--jet space.
 
\medskip

{\bf 1.10. A summary.} Deitmar's definition of the category of schemes
 over $F_1$ is, as he himself stresses 
in the opening paragraph of [De2],   a minimalistic one. 
It is quite transparent, but obviously does not allow one to treat some more sophisticated 
situations,
such as Soul\'e's scheme 1.7.2. In fact, the Theorem 4.1 of [De2] shows
that if $X$ is a connected integral $F_1$--scheme of finite type, then 
its lift to $\bold{C}$,  $X_{\bold{C}}$ consists of a finite union of 
mutually isomorphic toric varieties.

\smallskip

The richness of To\"en and Vaqui\'e's definition becomes apparent, when it is applied to
other basic symmetric monoidal categories. Especially remarkable is the
extension of $F_1$--geometry $\bold{S}_1$--$Sch$  which is the category of schemes relative to
the category  $(SEns, \times, *)$ of simplicial sets with direct product.
There is a canonical functor ``base extension''    $\bold{S}_1$--$Sch \to F_1$-$Sch$,
so that this geometry lies ``below''  $F_1$--geometry, in the same sense as
$F_1$--geometry lies below $\bold{Z}$--geometry. Another extension with great future
is the algebraic geometry over ``brave new rings''.

\smallskip

One outstanding problem is to extend cyclotomy to the homotopical framework.

\smallskip

This is an appropriate place to stress that in a wider context of [TV], or
eventually in noncommutative $F_1$--geometry, the spectrum of $F_1$ 
loses its privileged position of a final object of a geometric category.
For example, in noncommutative geometry, or in an appropriate category of stacks,
the quotient of this spectrum modulo the trivial action of a group
must lie below this spectrum.

\smallskip

Soul\'e's algebras $\Cal{A}_X$  are a very important element of the structure, 
in particular, because they form a bridge
to Arakelov geometry. Soul\'e uses concrete choices of  them in order to
produce  ``just right'' supply of morphisms, without  {\it a priori} constraining these choices formally.

\smallskip

However, these algebras appear as an {\it ad hoc} and somewhat
arbitrary supplement to the natural $F_1$--algebraic objects. Perhaps, 
a way to think about them is to imagine  {\it a possible
definition of  1--adic numbers}. 

\smallskip
Borger's context might lead to a progress in this direction.

\bigskip

\centerline{\bf  2.  Habiro's analytic functions of many variables:}

\smallskip

\centerline{\bf statements of results}

\medskip

{\bf 2.1. Notations.}  {\it Rings} in this and the next sections are associative, commutative and unital, unless a
context suggests otherwise. {\it Ring homomorphisms } are unital.
Letters $R,R_0, R_1 \dots$ denote rings, $q, q_0, q_1 ...$ are independent commuting variables.

\smallskip

Let $R$ be a ring, $\Cal{I}=\{I_{\alpha}\}$ a  family of ideals filtered by inclusion. The ring projective limit
$\varprojlim_{\alpha}R/I_{\alpha}$ is called the completion of $R$ with respect to $\Cal{I}$
and denoted $\wR_{\Cal{I}}$ or some version of this notation. When $\Cal{I}$ is (cofinal to)
the family of powers of one ideal $I$, the respective limit is called the $I$--adic completion.

\smallskip

We say that $R$ is $\Cal{I}$--  (resp. $I$--adically) {\it separated}, if $\cap_{\alpha}I_{\alpha}=\emptyset$.
Equivalently, the canonical homomorphism $R\to \wR_{\Cal{I}}$ is injective.
Example: $R=\bold{Z}$, $\Cal{I}$ any infinite filtering system.

\smallskip

When $q$ is considered as a  ``quantization parameter'', our quantized (Gaussian) versions of integers and
factorials are, as in [Ha2],
$$
\{N\}_q:= q^N-1,\ \{N\}_q!:= \{N\}_q\{N-1\}_q\dots \{1\}_q.
\eqno(2.1)
$$

\smallskip

Fix an integral domain $R_0$ of characteristic zero and put $R_n:= R_0[q_0,\dots, q_n],$
with natural embeddings $R_0\subset R_1\subset R_2\subset \dots $

\smallskip

Denote by $I_{n,N}\subset R_n$ the ideal  $(\{N\}_{q_1}!, \dots , \{N\}_{q_n}!)$, $N\ge 1.$
Clearly, $I_{n,N}\subset I_{n, N+1}$ so that the rings $R_n^{(N)}:=R_n/I_{n,N}$, $n\ge 1$ being fixed,
form an inverse system.

\medskip

{\bf 2.2.  Definition.} {\it The ring of Habiro's analytic functions of $n$ variables over $R_0$ is defined as} 
$$
\widehat{R}_n:= \varprojlim_N \, R_n^{(N)}.
$$

\medskip

{\bf  2.3. Taylor series of analytic functions.} Choose a vector of roots of unity
$\zeta =(\zeta_1,\dots, \zeta_n)$ such that all $\zeta_i$ are in $R_0$.
For any integer $M>0$, there exists $N_0=N_0(\zeta, M)$ such that 
$I_{n,N}\subset (q_1-\zeta_1,\dots ,q_n-\zeta_n)^M$ for all $N\ge N_0$.
In fact, $\{N\}_{q_i}!$ is divisible by any fixed monomial $(q_i-\zeta )^M, \zeta\in \mu$,
if $N$ is large enough.

\smallskip

The completion $\varprojlim_M R_n/ (q_1-\zeta_1,\dots ,q_n-\zeta_n )^M$
is $R[[q_1-\zeta_1,\dots ,q_n-\zeta_n ]].$

\smallskip

Therefore we obtain a ring homomorphism ``Taylor expansion at the point $\zeta$'':
$$
T_n(\zeta ):\,\wR_n \to R_0[[q_1-\zeta_1,\dots ,q_n-\zeta_n ]].
$$

\medskip

{\bf  2.3.1. Theorem.} {\it If $R_0$ is an integral domain, $p$--adically separated for all primes $p$,
then the same is true for $\wR_n$, and the Habiro--Taylor homomorphism $T_n(\zeta )$ is injective.}

\smallskip

More generally, let $F=\{F_1,\dots ,F_n\} \in \bold{Z} [q]$ be a family of monic polynomials in $R_0[q]$ whose all roots are
roots of unity. Denote by  $(F)$ the ideal generated $F_1(q_1),\dots ,F_n(q_n)$ in $R_n.$
In place of the formal series ring above, we can consider the completion
$$
\wR_F:= \varprojlim_M R_n/(F)^M
$$
and the respective Taylor expansion homomorphism:
$$
T_n(F):\, \wR_n \to \wR_F.
$$

\medskip

{\bf  2.3.2. Theorem.} {\it If $R_0$ is an integral domain, $p$--adically separated for all $p$,
$R[[F]]$ is as well $p$--adically separated, and the homomorphism $T_n(F)$ is injective.}

\smallskip

K.~Habiro proved these results, as well as their generalizations, for $n=1$, and we build upon
his proof.

\medskip

{\bf 2.3.3. Differential calculus.} Divided powers of partial derivatives with respect to $q_k$ are continuous in
the linear topologies  generated by  $I_{n,N}$, resp. by all    all  $(q_1-\zeta_1,\dots ,q_n-\zeta_n)^M$.
Hence these derivatives make sense in  $\widehat{R}_n$, and their
values at $(\zeta_1, \dots ,\zeta_n)$ are the Taylor coefficients 
of the respective series.  

\smallskip

(In order to check the continuity with respect to  $I_{n,N}$ it suffices to notice
that as $N$ tends to infinity, $\{N\}_q!$ as a polynomial of $q$ vanishes
at a growing set of roots of unity with infinitely growing multiplicity
at each root of unity. Taking a derivative of such a sequence of polynomials does
not destroy this property).

\smallskip

Thus we can develop for $\wR_n$ the conventional formalism
of tangent and cotangent modules, differential forms etc.

\medskip

{\bf 2.4. Elements of $\wR_n$ as functions on roots of unity.}  Let $R^{\prime}_0\supset R_0$
be an integral domain flat over $R_0$ and containing all roots of unity (that is, all cyclotomic
polynomials  $q^n-1$ completely split in $R^{\prime}_0$). Denote by $\mu$ the set of all roots
of unity in $R^{\prime}_0.$ Choose $\zeta:=(\zeta_1,\dots ,\zeta_n)\in \mu^n.$
Any element of $R_n$, being a polynomial in $(q_1,\dots ,q_n)$, takes a certain value at $\zeta$ belonging to $R_0^{\prime}$.
If $N\ge N_0(\zeta )$, all elements of $I_{n,N}$ vanish at $\zeta$.
Hence any element $f\in \wR_n$ defines a map $\bar{f}:\, \mu^n\to R_0^{\prime}$.
This map is $R_0$--linear and  compatible with pointwise addition and multiplication of functions.
\smallskip

Besides assuming that $R_0$ is $p$--adically separated for all primes $p$,
impose the following separatedness condition:
{\it for any infinite sequence of pairwise distinct primes $p_1, \dots , p_k, \dots $,
we have}
$$
\cap_{m=1}^{\infty} Rp_1\dots p_m =\{0\}.
\eqno(2.2)
$$

\medskip

{\bf 2.4.1. Theorem.} {\it  Under  these assumptions, the map $f\mapsto \bar{f}$ is injective.}

\smallskip

One can also formulate this statement without adjoining to $R_0$ roots of unity.

\medskip

{\bf 2.4.2. Theorem.} {\it 
The natural map $\wR_n\to \prod_{m=1}^{\infty}\wR_n\,\roman{mod}\, (\Phi_m(q_1),\dots ,\Phi_m(q_n))$    is injective.}

\smallskip

For $n=1$, these results were established by K.~Habiro. He has also shown that vanishing
of $\bar{f}$ on certain sufficiently large subsets of $\mu$ suffices to establish the
vanishing of $f$. 

\smallskip

More precisely, {\it Habiro's topology} on the set $\mu$ of all roots of unity is defined as follows
(cf. [Ha2], 1.2).

\smallskip

Two roots of unity $\xi,\eta$ are called {\it adjacent}, if  $\xi\eta^{-1}$
is of order $p^m$, $m\in\bold{Z}$, $p$ a prime; or equivalently, if $\xi -\eta$ is not a unit
(as an algebraic number). Clearly, the action of $Gal\,(\overline{\Q}/\Q )$ 
preserves adjacency.

\medskip

{\bf 2.4.3. Definition.} {\it A subset $U\subset \mu$ is called open, if
for any point $\xi\in U$, all except of finitely many $\eta\in \mu$ , adjacent
to $\xi$, belong to $U$.}

\smallskip

The Galois action is continuous in this topology, in marked contrast to the topology induced from $\bold{C}.$

\smallskip

Let now $\mu^{\prime}$ be an infinite set of roots of unity. A point
$\xi\in \mu^{\prime}$ is a limit point of $\mu^{\prime}$,
if for any open neighborhood $U$ of $\xi$ we have $\mu^{\prime}\cap (U\setminus \xi)\ne \emptyset.$
In Habiro's topology, this means that $\mu^{\prime}$ contains infinitely many points, adjacent to $\xi$.

\medskip

{\bf 2.4.4. Theorem.} {\it Under the notations and assumptions of Theorem 2.3.1, 
let $\nu =\nu_1\times\dots\times \nu_n \subset \mu^n$ be a set, such that each $\nu_i\subset \mu$ 
has a limit point.  Let $f\in \wR_n$.
If the restriction $\bar{f} |_{\nu}$ is identical zero, then $f=0$.}

\smallskip

In the next section, we will prove this last result; Theorems 2.4.1 and 2.4.2
follow from it.

\medskip

{\bf 2.5. Analogs of Habiro's functions on the arithmetic axis and analytic
continuation.}
The Habiro ring of one variable $\varprojlim_N \bold{Z}[q]/(\{N\}_q!)$
``is'' the lift to $\bold{Z}$ of an imaginary ring $\varprojlim_N F_1[q]/(\{N\}_q!)$.

\smallskip

Along the arithmetical axis, the straightforward analog of the latter exists: this is the topological ring
of profinite integers $\widehat{\bold{Z}}:=\varprojlim_N \bold{Z}/(N!)$.
Its elements can be uniquely represented by  infinite series
$\sum_{n=1}^{\infty} c_n n!$ where $c_n$ are integers with $0\le c_n\le n,$
cf. [van D].

\smallskip

H.~Lenstra in [Le] discusses profinite Fibonacci numbers:  continuous extrapolation
to $n\in\widehat{\bold{Z}}$ of the Fibonacci function $n\mapsto u_n$.

\smallskip

An analog of the profinite number $1+\sum_{n=1}^{\infty} (-1)^nn!$ is the remarkable
example of Habiro function of one variable 
$$
1+\sum_{n=1}^{\infty} (-1)^n\{n\}_q! = 1 +\sum_{n=1}^{\infty} (1-q)\dots (1-q^n).
$$
As a function on roots of unity, it  emerged in a work of M.~Kontsevich on Feynman integrals 
(talk at MPIM, 1997). Don Zagier in [Za] proved that its values, as well 
as values of its derivatives,  are radial limits of the function (resp. its derivatives)
holomorphic in the unit circle
$$
\dfrac{1}{2}\sum_{n=1}^{\infty} n\chi (n)q^{(n^2-1)/24},
$$
where $\chi$ is the quadratic character of conductor 12.

\medskip

{\bf 2.6. Habiro's functions on $F_1$--schemes.} Let $X$ be an
$F_1$--scheme in the sense of one of the definitions from sec. 1.
Let $(x_1,\dots ,x_n)$ be a finite family of local cyclotomic
coordinates on $X$. For any ring $R$ as in 1.7. (ii),
denote by $U(R)\subset X(R)$ the set
of cyclotomic points, at which all $x_i$ are defined and take 
non--zero values.

\smallskip

Consider an analytic function
$f\in \wR_n$ in the sense of Habiro. This function
then defines a map
$$
f_R:\, U(R)\to R,\, f_R(r):= \bar{f}(x_1(r),\dots ,x_n(r)),
$$
with evident functorial properties.  

\smallskip

In an appropriate setting such functions must be local sections
of a global sheaf. I hope to return to this problem in another
paper. Here I will restrict myself to the following
observations.

\smallskip

(i) We have to exclude zero values, because  $q_1$ is invertible in $\wR_1$,
and hence each monomial $q_1^{m_1}\dots q_n^{m_n}$ is invertible
in $\wR_n$. In fact, 
$$
q^{-1}=1+\sum_{n=1}^{\infty} (-1)^n q^n\{n\}_q!,
$$
see [Ha1], Proposition 7.1.

\smallskip

(ii) From the perspective of this paper, it seems quite
natural to consider localizations with respect to functions
such as $q_1^{m_1}\dots q_n^{m_n}-1$, deleting sets of roots
of unity closed in Habiro's topology. However, such functions
are generally not cyclotomic coordinates. This runs counter the spirit
of To\"en--Vaqu\'e's definitions, and requires rethinking of their
framework.

\bigskip

\centerline{\bf  3.  Habiro's analytic functions of many variables:}

\smallskip

\centerline{\bf proofs and generalizations}

\medskip

{\bf 3.1. The case $n=1$.} Assuming that a ring
$R$ is $\Cal{I}$--separated for each member $\Cal{I}$ of some set of filters $\Cal{S}_R$,
we can deduce that the ring $R[q]$, and certain its completions, are separated
with respect to the members of another set of filters, say $\Cal{S}_{R[T]}$.
Results of this type are collected and proved in [Ha1].
They will allow us to perform inductive steps, passing from $n$ to $n+1$,

\smallskip

{\bf 3.2. Proof of the Theorem 2.3.1.} We will perform induction on $n$, using Habiro's theorem
for $n=1$ ([Ha1], Theorem 5.2) as the basis of induction. 

\smallskip

Assuming the theorem proved for $\wR_n$, we will proceed by decomposing the Taylor series 
map  $\wR_{n+1} \to R_0[[q_1-\zeta_1,\dots ,q_n-\zeta_n,q_{n+1}-\zeta_{n+1}]]$ into the product of two ring
homomorphisms and checking injectivity of each one:
$$
\wR_{n+1}\overset\alpha\to{\to} \wR_n[[q_{n+1}-\zeta_{n+1}]]  \overset\beta\to{\to} R_0[[q_1-\zeta_1,\dots , q_n-\zeta_n]]\, [[q_{n+1}-z_{n+1}]] 
$$

\smallskip 

Now we will define the arrows $\alpha ,\beta$ and check their properties.

\smallskip

The arrow $\beta$ is continuous in $(q_{n+1}-\zeta_{n+1})$--adic topology, acts identically on $q_{n+1}-\zeta_{n+1}$, 
and sends each element of $\wR_n$ to its Taylor series at $(\zeta_1,\dots ,\zeta_n).$ In view of the inductive
assumption, $\beta$ is injective.
  \smallskip
  
  To define $\alpha$,
  consider an element $g\in \wR_{n+1}$. It can be represented as
  the limit of a sequence of polynomials $g_1, g_2, \dots ,g_N, \dots$, where $g_i\in R_0 [q_1,\dots ,q_{n+1}]$
  such that $g_{N+1} \equiv  g_N \,\roman{mod}\, I_{n+1,N}$. 
  
  \smallskip
  
  From the definition it follows that
  $$  
  I_{n+1,N}= I_{n,N} [q_{n+1}] + R_{n+1}\cdot\{N\}_{q_{n+1}}!
  $$
 Therefore,
 $$
 g_{N+1}=g_N+ i_N + r_N\cdot \{N\}_{q_{n+1}}!, 
 \eqno(3.1)
 $$
 where
 $$
  i_N\in I_{n,N} [q_{n+1}], \ r_N  \in R_{n+1}.
 $$
 
Now consider a point $(\zeta_1, \dots ,\zeta_{n+1})$ as above.
Clearly,
$$
I_{n,N} [q_{n+1}]=I_{n,N} [q_{n+1}-\zeta_{n+1}].
$$

\smallskip
 Write $g_N, i_N, \{N\}_{q_{n+1}}!$ as  polynomials
in $q_{n+1}-\zeta_{n+1}$ with coefficients in $R_n$.  When $N$ becomes large enough,  $\{N\}_{q_{n+1}}!$ starts with arbitrary large
power of $q_{n+1}-\zeta_{n+1}$. Therefore for any given $M$, the coefficient at $(q_{n+1}-\zeta_{n+1})^M$
in  $g_{N+1}$ is the same as in $g_N+i_N$ if $N\ge N_1(M,\zeta ).$ Hence the sequence
of these coefficients ($M$ being fixed and $N$ growing) converges to a certain element $a_M\in \wR_n$. 

\smallskip

Put $\alpha (g) := \sum_{M=0}^{\infty} a_M(q_{n+1}-\zeta_{n+1})^M$.
One can routinely check that $\alpha (g)$ depends only on $g\in\wR_{n+1}$ and not on the system $(g_N)$
chosen to represent $g$. Moreover, we get a ring homomorphism
$$
\alpha :\, \wR_{n+1}\to \wR_n [[q_{n+1}-\zeta_{n+1}]] .
\eqno(3.2)
$$
\smallskip

Let us check that $\alpha$ is injective. In fact, take a nonzero element $g=\lim g_N$.
Then there exist arbitrarily large $N$ such that $g_N\notin I_{n+1,N}.$
Representing $g_N$ as a polynomial in $q_{n+1}-\zeta_{n+1}$ with coefficients in $R_n$, we can find in this polynomial a 
coefficient, not belonging to $I_{n,N}$. In the limit, it will produce a nonvanishing $a_M$.

\smallskip

Finally, $\beta\circ \alpha =T_{n+1}( \zeta )$ by construction.

\smallskip

\medskip

{\bf 3.3.  Proof of the Theorem 2.4.4.} We first remark, that the case $n=1$ is essentially covered 
by the Theorem 6.1 of [Ha1],
if one weakens the assumption $R_0\subset \bar{\Q}$ in the statement of this Theorem. In fact, this 
assumption is used only at the end of the proof, in order to ensure the validity of the separatedness
condition (2.2). Instead, we will simply postulate (2.2)  for $R_0$,
and then deduce it for each $\wR_n$ using the Taylor embedding of $R_n$ into
$R_0[[q_1-\zeta_1,\dots ,q_{n+1}-\zeta_{n+1}]].$

\smallskip

To pass from $n$ to $n+1$,  I will start with the following remarks.

\smallskip

Let $R$ be a ring endowed with a filtering family of ideals $\Cal{I} = \{I_{\alpha}\}.$
Consider the following two families of ideals in the polynomial ring $R[q]$:

\smallskip

(i) $\Cal{I}_1 :=\{ I_{\alpha}[q]+ (\{N\}_q!)\,|\, \alpha , N\ \roman{arbitrary} \}.$

\smallskip

(ii) $\Cal{I}_2 :=\{ (\{N\}_q!)\,|\,  N\ \roman{arbitrary} \}.$

\smallskip

Denote by $R[q]^{\,\widetilde{}}$ (resp. $R[q]^{\,\widehat{}}$) the completion of $R[q]$ with respect to
$\Cal{I}_1$ (resp. $\Cal{I}_2$.

\smallskip

For any $N$ and $\alpha$, we have natural surjections
$$
R[q]/(\{N\}_q!)\to R[q]/( I_{\alpha}[q]+ (\{N\}_q!)).
$$
Passing to the limit, we get a canonical surjection
$$
\varphi :\, R[q]^{\,\widehat{}} \to R[q]^{\,\widetilde{}}.
$$

\smallskip

{\bf 3.3.1. Lemma.}  {\it Consider the case $R=\wR_n$, $\Cal{I} = \{\widehat{I}_{n,N})$
where
$$
\widehat{I}_{n,N} := (\{N\}_{q_1}!,\dots , \{N\}_{q_n}!)\subset \wR_n .
$$
Then the homomorphism
$$ 
\varphi :\,\wR_n[q_{n+1}]^{\,\widehat{}}\to \wR_n[q_{n+1}]^{\,\widetilde{}} =\wR_{n+1}
$$
is an isomorphism.}

\smallskip

{\bf Proof.} It suffices to check that $Ker\,\varphi = \{0\}.$  In fact, as in 3.2, we have an injection
$$
\alpha :\,\wR_n[q_{n+1}]^{\,\widetilde{}}\to \wR_n[[q_{n+1}-1]]
$$
and an one--variable Taylor series injection
$$
T:\, \wR_n[q_{n+1}]^{\,\widehat{}}\to\wR_n[[q_{n+1}-1]].
$$
By construction, $\alpha\circ\varphi =T$, hence $\varphi $ is an injection as well.

\medskip

{\bf 3.3.2. End of proof of the Theorem 2.4.4.} Suppose now that $g\in \wR_{n+1}$ vanishes
at all points $(\zeta_1,\dots ,\zeta_{n+1})$, $\zeta_i\in\nu_i\subset \mu$,
each of $\nu_i$ having a limit point. To simplify notation, assume that all roots
of unity are in $R_0$.

\smallskip

The evaluation of $g$ at $(\zeta_1,\dots ,\zeta_{n+1})$ can be decomposed into the composition
of two arrows:
$$
ev_{(\zeta_1,\dots ,\zeta_n)}\circ ev_{\zeta_{n+1}} :\,      \wR_{n+1}\to \wR_n\to R_0 ,
$$
where the first arrow $ev_{\zeta_{n+1}}$ is obtained by taking the constant term in $\alpha (g)$, (3.2), and the second one
is the evaluation at $(\zeta_1 ,\dots ,\zeta_n)$.

\smallskip

First, fix $(\zeta_1,\dots ,\zeta_n)$ and vary $\zeta_{n+1}\in \nu_{n+1}.$
We have already identified $\wR_{n+1}$ with $\wR_n[q_{n+1}]^{\,\widehat{}}$
in a way which is clearly compatible with evaluation maps.

\smallskip

From the Habiro Theorem 6.1, [Ha1], we obtain that 
$$
ev_{(\zeta_1,\dots ,\zeta_n)}(g)=0
$$
for all 
$$
(\zeta_1,\dots ,\zeta_n)\in \nu_1\times\dots\times \nu_n.
$$
By the inductive assumption, $g=0$. This finishes the proof.

\medskip

{\bf 3.4. General monoids, coordinate independence, and functorality.}  Let $M$ be a commutative monoid with unit.

\smallskip

We can consider the completion $R_0^{\prime}[M]$  of $R_0 [M]$ with respect to the system of ideals
$I_N$, where $I_N$ is generated by all elements $\{N\}_m!:=(m^n-1)\dots (m-1)$ for $m\in M$.

\smallskip

Obviously, any morphism $\psi :\, M\to N$ induces the respective morphism
of the completed rings. In particular, the diagonal morphism $M\to M\times M$
produces a structure of Hopf algebra on $R_0[M]$ and its completed version
on $R_0^{\prime}[M]$.

\smallskip

As K.~Habiro noticed in a message to the author (Aug. 23, 2008), applying this construction to $M=\bold{Z}^n$,
we get precisely $\wR_n$ (if $q_i$ corresponds to the basic vector $(0,\dots ,1, 0, \dots , 0)$,
with 1 at $i$--th place.) Since $q_i$ are invertible in $\wR_n$, we could as well start
with $R_0[q_1, q_1^{-1},\dots ,q_n,q_n^{-1}]$, but it seemed  more natural to me to deduce the invertibility
at the end of the construction.

\bigskip

\centerline{\bf 4. Schemes with natural cyclotomic coordinates:}

\smallskip
\centerline{\bf Witt vectors and moduli spaces}

\medskip

In this section we treat two disjoint constructions.

\smallskip

{\bf 4.1. Witt functors.} The (big) Witt ring scheme $W$ can be defined
as an infinite dimensional affine space $Spec\,\bold{Z} [u_1,u_2,u_3,\dots ]$,
whose polynomial algebra of functions $A$ is endowed 
with two homomorphisms $A\to A\otimes A$, ``coaddition'' $\bold{\alpha}$ and ``comultiplication''
$\bold{\mu}$.

\smallskip

The functor of its $R$--points, for a variable commutative ring $R$,
set theoretically is $W(R) = \prod_{k=1}^{\infty} R$ where the $k$--th coordinate
of the product is the value of $u_k$ at the respective $R$--point.
The maps $\bold{\alpha}$ and $\bold{\mu}$ induce on $W(R)$ the structure of 
commutative ring, functorial in $R$. This structure can be described quite
explicitly, if we use in place of $\{u_k\}$ the  ``ghost coordinates''
$$
q_n:= \sum_{d|n} du_d^{n/d}.
$$
In this coordinates, $\bold{\alpha}$ and $\bold{\mu}$ induce respectively
componentwise addition and multiplication (cf. [Haz], sec. 9 and 14, in particular (14.3)).
 
 \smallskip
 
 The $N$--truncated Witt scheme $W^{(N)}$ is obtained if we apply this to the subring
  $\bold{Z}[u_1,\dots ,u_N]$ with induced $\bold{\alpha}_N$ 
 and $\bold{\mu}_N$. For a prime $p$, the scheme $W_p$ is obtained by
 taking the subring generated by all $u_{p^k}$, $k\ge 0.$ The
 truncated version $W_p^{(N)}$ jumps only at powers of $p$ as well.
In this way we get quotient
functors of the Witt functor, valued
in commutative algebras.
 
 \smallskip
 
 In place of subrings, one can consider  quotients by the  ideals generated by the
 complementary coordinates.

\smallskip

{\bf 4.1.1. Definition. }{\it The (truncated) Witt gadget $\Cal{W}^{(N)}$ is
defined by the following data:}

\smallskip

(i) $\Cal{W}^{(N)}_{\bold{Z}}:= W^{(N)}. $

\smallskip

(ii)   {\it For a ring $R$ as in 1.7, (ii), the subfunctor of cyclotomic points
$\Cal{W}^{(N)}_{\bold{Z}}(R)$ of $W^{(N)}(R)$ is defined as consisting
of points, whose ghost coordinates are 0 or roots of unity. }

\smallskip

Thus, ghost coordinates are cyclotomic coordinates in the sense of 1.7.1.

\medskip

{\bf 4.2. Moduli spaces $\overline{L}_{0;2,B}$.}  In an ideal world,  not only
schemes allowing ``finite combinatorial'' description  ([So], p. 217) must be extensions of objects
over $F_1$, but perhaps  ``all''  rigid structures as well. An obvious challenge
is presented by $\overline{M}_{0,n}$, moduli spaces of stable curves of genus
zero with $n$ marked points forming the basic operad of quantum cohomology. 

\smallskip

As the first approximation, we look in this subsection to some moduli spaces introduced in 
[LoMa1] and studied further in [LoMa2] and [Ma2]. Generally, they parametrize
curves of genus $g$, with marked points, a part of which
(carrying ``black''  labels) being allowed to merge between them, although
not with singular or ``white labeled'' points. There is an appropriate notion
of stability and  a representability theorem. (Both were vastly generalized in the study 
[BaMa].)

\smallskip

Here we will focus on  the case of genus zero, two white points and arbitrary ($\ge 1$)
number of black points. The resulting moduli spaces turn out to be toric,
based upon permutohedral fans. Therefore they are certainly
lifts to $\bold{Z}$ of toric $F_1$--schemes. We discuss which of the canonical
morphisms between them  
descend to $F_1$.

\smallskip

It is convenient to label the black points by elements of a finite set $B$,
carrying no additional structure (rather than, say, by $ \{1,\dots ,n\}$,
which suggests a complete order on labels).

\smallskip

Below we will give a toric description of the respective moduli space that we will now denote
$\overline{L}_B.$ For proofs, see  [LoMa1].

\medskip

{\bf 4.3. Partitions.} A partition $\{\sigma\}$ of  a finite
set $B$ is {\it a totally
ordered set of non--empty subsets of $B$ whose union is
$B$ and whose pairwise intersections are empty.}
If a partition consists of $N$ subsets, it is called
$N$--partition. If its components are denoted $\sigma_1,\dots ,\sigma_N$,
or otherwise listed, this means that they are listed
in their structure order.

\smallskip

Let $\tau$ be an $N+1$--partition of $B$. If $N\ge 1,$ it determines
a well ordered family of $N$ 2--partitions $\sigma^{(a)}$:  
$$
\sigma^{(a)}_1:=\tau_1\cup\dots\cup\tau_{a},\
\sigma^{(a)}_2:=\tau_{a+1}\cup\dots\cup\tau_{N+1},\ a=1,\dots ,N\, .
\eqno(4.1)
$$
In reverse direction, call a family of 2--partitions $(\sigma^{(i)})$ 
 {\it good} 
if for any $i\ne j$ we have $\sigma^{(i)}\ne \sigma^{(j)}$
and either  $\sigma^{(i)}_1\subset \sigma^{(j)}_1,$
or $\sigma^{(j)}_1\subset \sigma^{(i)}_1.$ Any good family
is naturally well--ordered by the relation 
$\sigma^{(i)}_1\subset \sigma^{(j)}_1$, and we will
consider this ordering as a part of the structure. If a good family of 2--partitions
consists of $N$ members, we will usually choose
superscripts $1,\dots ,N$ to number these partitions
in such a way that $\sigma^{(i)}_1\subset \sigma^{(j)}_1$
for $i<j.$

\smallskip

Such a good family produces one $(N+1)$--partition $\tau$:
$$
\tau_1:=\sigma_1^{(1)},\ \tau_2:=\sigma_1^{(2)}\setminus
\sigma_1^{(1)},\ \dots ,\  
\tau_N:=\sigma_1^{(N)}\setminus
\sigma_1^{(N-1)},\ \tau_{N+1}=\sigma_2^{(N)}.
\eqno(4.2)
$$ 
This correspondence between good $N$--element families of
2--partitions and $(N+1)$--partitions is one--to--one, because
clearly $\sigma_1^{(i)}=\tau_1\cup\dots\cup\tau_i$ for $1\le i\le N.$

\medskip

{\bf 4.4. The fan $F_B$.} Now we will describe
a fan $F_B$ in the space $N_B\otimes{\bold{R}}$, where
$N_B:=\roman{Hom}\,(\bold{G}_m,T_B)$, $T_B:=\bold{G}_m^B/\bold{G}_m$.
Clearly, $N_B$ can be canonically identified with
$\bold{Z}^B/\bold{Z}$, the latter subgroup being embedded
diagonally. Similarly, $N_B\otimes{\bold{R}}=\bold{R}^B/\bold{R}$.
We will write the vectors of this space (resp. lattice) as functions
$B\to \bold{R}$ (resp. $B\to \bold{Z}$) considered modulo
constant functions. For a subset $\beta\subset B$, let $\chi_{\beta}$
be the function equal 1 on $\beta$ and 0 elsewhere.

\medskip

{\bf 4.4.1. Definition.}  {\it The fan $F_B$ consists
of the following $l$--dimensional cones $C(\tau )$ labeled by 
$(l+1)$--partitions $\tau$ of $B$.

\smallskip

If $\tau$ is the trivial 1--partition, $C(\tau )=\{0\}$.

\smallskip

If $\sigma$ is a 2--partition, $C(\sigma )$ is generated
by $\chi_{\sigma_1}$, or, equivalently, $-\chi_{\sigma_2}$,
modulo constants.

\smallskip

Generally, let $\tau$ be an $(l+1)$--partition, and 
$\sigma^{(i)},\,i=1,\dots,l$,
the respective good family of 2--partitions (4.1). Then
$C(\tau )$ as a cone is generated by all $C(\sigma^{(i)})$. }

\medskip

{\bf  4.5. Toric varieties  $\overline{L}_B$ and forgetful morphisms.} We denote by
$\overline{L}_B$ the variety associated with the fan $F_B$.
It is smooth and proper, in fact projective.

\smallskip

Assume that $B\subset B^{\prime}$. Then we have the projection
morphism $\bold{Z}^{B^{\prime}}\to\bold{Z}^B$ which induces the morphism
$f^{B^{\prime},B}:\,N_{B^{\prime}}\to N_B.$ It satisfies the following property:
 for each cone 
$C(\tau^{\prime})\in F_{B^{\prime}}$,
there exists a cone $C(\tau )\in F_{B}$  such that
$f^{B^{\prime},B}(C(\tau^{\prime}))\subset C(\tau ).$ In fact,
$\tau$ is obtained from $\tau^{\prime}$ by deleting
elements of $B^{\prime}\setminus B$ and then deleting the
empty subsets of the resulting partition of $B$.

\smallskip

Therefore, we have a morphism $f^{B^{\prime},B}_*:\,  
\overline{\Cal{L}}_{B^{\prime}}\to\overline{\Cal{L}}_B$
 which we will call {\it forgetful one} (it forgets
elements of $B^{\prime}\setminus B$).The forgetful morphism is flat, because locally in toric
coordinates it is described as adjoining  variables and localization.

\medskip

{\bf 4.6. $\overline{L}_B$ as families of curves with two white and $B$ black points.}
This structure can be defined in terms of forgetful morphisms forgetting just one point  $B$.
Let $B\subset B^{\prime}$, $\roman{card}\,B^{\prime}\setminus B =1.$
\smallskip

We start with describing structure sections.

\smallskip

In order to define the two white sections of the forgetful
morphism, consider two partitions $(B^{\prime}\setminus B,B)$
and $(B,B^{\prime}\setminus B)$ of $B^{\prime}$ and the respective
closed strata. The forgetful morphism
restricted to these strata identifies them with $\overline{\Cal{L}}_B$.
We will call them $x_0$ and $x_{\infty}$ respectively.

\smallskip

Finally, to define the $j$--th black section, $j\in B$,
consider the morphism of lattices $s_j:\,N_B\to N_{B^{\prime}}$
which extends a function $\chi$ on $B$ to the function
$s_j(\chi )$ on $B^{\prime}$ taking the value $\chi (j)$ at the forgotten
point. This morphism satisfies the following condition:
each cone $C(\tau )$ from $F_B$ lands in an appropriate cone 
$C(\tau^{\prime})$ from $F_{B^{\prime}}$.
Hence we have the induced morphisms $s_{j*}:\,\overline{\Cal{L}}_B\to
\overline{\Cal{L}}_{B^{\prime}}$ which obviously are sections.
Moreover, they do not intersect $x_0$ and $x_{\infty}$.

\smallskip

{\bf 4.6.1. Proposition.} {\it With the notations and assumptions above, the forgetful morphism
is a universal family of (painted stable) marked curves of genus zero with
two white points and $B$ black points.}

\smallskip

In order to see the structure of fibers of the forgetful morphism,
one should notice that the inverse image of any point $x\in\Cal{L}_{\tau}$ is acted
upon by the multiplicative group $\bold{G}_m =\roman{Ker}\,(T_{B^{\prime}}\to T_B)$.
This action breaks the fiber into a finite number of orbits
which coincide with the intersections of this fiber
with various $\Cal{L}_{\tau^{\prime}}$ described above.
When $\tau^{\prime}$ is obtained by adding the forgotten
point to one of the parts, this intersection is a torsor
over the kernel, otherwise it is a point. As a result, we get
that the fiber is a chain of $\bold{P}^1$'s, whose components are labeled
by the components of $\tau$ and singular points by the neighboring
pairs of components.

\medskip

{\bf 4.7. Clutching morphisms.} They are morphisms of the type
$\overline{L}_{B_1}\times\overline{L}_{B_2}\to\overline{L}_{B_1\coprod B_2}$
whose fiberwise description is this: glue $\infty$ of the first curve to
$0$ of the second curve. They admit an obvious toric description.

\smallskip

About their operadic role, see [Ma2].

\smallskip

{\bf 4.8. Proposition.} {\it  Forgetful and clutching morphisms descend
to the $F_1$--models of the toric varieties $\overline{L}_B$.}

\bigskip

\centerline{\bf References}

\smallskip
 
[BaMa] A.~Bayer, Yu.~Manin. {\it Stability conditions, wall--crossing and weighted 
Gromov--Witten invariants.}
e-print math.AG/0607580

\smallskip

[Bor1] J.~Borger. {\it The basic geometry of Witt vectors.}  e-print math.AG/0801.1691

\smallskip

[Bor2] J.~Borger. {\it $\Lambda$--rings and the field with one element. Preliminary version.} Preprint, 2009.

\smallskip

[BorS] J.~Borger, B.~de Smit. {\it Galois theory and integral models
of $\Lambda$--rings.} Bull.~Lond.~Math.~Soc., 40 (2008), 439--446.
e--print arxiv:0801.2352
\smallskip

[BoCo] J.~Bost, A.~Connes. {\it Hecke algebras, type III factors and phase
transitions with spontaneous symmetry breaking in number theory.}
Selecta Math., 3 (1995), 411--457.

\smallskip

[CC] A.~Connes, C.~Consani. {\it On the notion of geometry over $F_1$.}
e--print arxiv:0809.2926

\smallskip
[CCMa1] A.~Connes, C.~Consani, M.~Marcolli. {\it Noncommutative geometry and 
motives: the thermodynamics of endomotives.}
Advances in Math., 214 (2) (2007), 761--831.

\smallskip

[CCMa2] A.~Connes, C.~Consani, M.~Marcolli. {\it Fun with $F_1$.}  e--print

math.AG/0806.2401

\smallskip

[CCMa3]  A.~Connes, C.~Consani, M.~Marcolli. {\it  The Weil proof and the geometry of the
ad\`eles class space.} To appear in ``Algebra, arithmetic and Geometry -- Manin Festschrift'', 
Progress in Math., Birkh\'auser, 2009. e-print math/0703392
\smallskip

[De1] A.~Deitmar. {\it Schemes over $F_1$.} In: Number Fields and Function Fields -- Two Parallel
Worlds.  Ed. by G.~ van der Geer, B.~Moonen, R.~Schoof.  Progr. in Math, vol. 239, 2005 . e-print math.NT/0404185

\smallskip
[De2]  A.~Deitmar. {\it $F_1$--schemes and toric varieties.}  e-print math.NT/0608179
\smallskip
[De3]  A.~Deitmar. {\it Remarks on zeta functions and $K$--theory over $F_1$.} Proc. Jap. Ac. Sci,
ser. A, vol. 82:8 (2007), 141--146.

\smallskip

[Du] N.~Durov. {\it New approach to Arakelov geometry.} e-print math/0704.2030

\smallskip

[Go] V.~Golyshev. {\it The canonical strip I.} Preprint, May 2008.
\smallskip
[Gr]  D.~Grayson. {\it $SK_1$ of an interesting principal ideal domain.} Journ. of Pure and Appl. Algebra,
20 (1981), 157--163.
\smallskip
[Ha1] K.~Habiro. {\it Cyclotomic completions of polynomial  rings.} Publ. RIMS, Kyoto Univ., 40 (2004), 1127--1146.
\smallskip
[Ha2] K.~Habiro. {\it A unified Witten--Reshetikhin--Turaev invariant for integral homology spheres.}
Inv. Math., 171 (2008), 1--81.
\smallskip

[Har] S.~Haran. {\it Non--additive geometry.} Comp. Math., 143 (2007), 618--688.

\smallskip

[Haz] M.~Hazewinkel. {\it Witt vectors. Part I.} e-print math/0804.3888

\smallskip

[Ka] M.~Kapranov. {\it Some conjectures on the absolute direct image.} Letter, May 1995.

\smallskip
[KaS]. M.~Kapranov, A.~Smirnov. {\it Cohomology determinants and reciprocity laws:
number field case.}  Unpublished manuscript, 15 pp.

\smallskip

[Ku] N.~Kurokawa. {\it Zeta functions over $F_1$.} Proc. Jap. Ac. Sci., Ser. A, vol. 81 (2005), 180--184.
\smallskip
[Law]  R.~J.~Lawrence. {\it Witten--Reshetikhin--Turaev invariants of 3--manifolds as holomorphic functions.}
In: Geometry and Physics (Aarhus, 1995), LN in Pure and Appl. Math., 184 (1997), Dekker, NY.

\smallskip

[LawZ] R.~J.~Lawrence, D.~Zagier. {\it Modular forms and quantum invariants of 3--manifolds.} Asian J. Math., 3 (1999),
93--107.

\smallskip

[Le] H.~Lenstra. {\it Profinite Fibonacci numbers.}  Nieuw Arch. Wiskd. (5)  6  (2005),  no. 4, 297--300.

\smallskip

[LoMa1]  A.~Losev, Y.~Manin. {\it New moduli spaces of pointed curves and pencils
of flat connections.}
  Fulton's Festschrift,
Michigan Journ. of Math., vol. 48, 2000, 443--472. e-print math.AG/0001003

\smallskip

[LoMa2]  A.~Losev, Y.~Manin. {\it Extended modular operad. }
In: Frobenius Manifolds,
ed. by C. Hertling and M. Marcolli,  Vieweg \& Sohn Verlag,
Wiesbaden, 2004, 181--211. e-print 
math.AG/0301003

\smallskip
[Ma1] Yu.~Manin. {\it Lectures on zeta functions and motives (according to Deninger and Kurokawa).}
Ast\'erisque 228:4 (1995), 121--163.

\smallskip

[Ma2]  Yu.~Manin. {\it Moduli stacks $\overline{L}_{g,S}.$}  Moscow Math. Journal,
4:1 (2004), 181--198. e--print math.AG/0206123

\smallskip

[Marc] M.~Marcolli. {\it Cyclotomy and endomotives.} e--print  arxiv:0901.3167

\smallskip

[Mart1] F.~Marty. {\it Relative Zariski open objects.} Preprint.

\smallskip

[Mart2] F.~Marty. {\it Smoothness in relative geometry.} Preprint.

\smallskip

[RV] F.~Rodriguez--Villegas. {\it  On the zeroes of certain polynomials.} Proc. AMS, 130:8 (2002), 2251--2254.
\smallskip
[Sm1]  A.~L.~Smirnov.  {\it Hurwitz inequalities for number fields. (Russian).}  Algebra i Analiz  4  (1992),  no. 2, 186--209; 
 translation in  St. Petersburg Math. J.  4  (1993),  no. 2, 357--375.
 
 \smallskip
 
 [Sm2]  A.~L.~Smirnov. {\it Absolute determinants and Hilbert symbols.} Preprint MPI 94/72, Bonn, 1994.
\smallskip
[So] C.~Soul\'e. {\it Les vari\'et\'es sur le corps \`a un \'el\'ement.} Mosc. Math. J. 4:1 (2004), 217--244.
\smallskip

[Ti] J.~Tits. {\it Sur les analogues alg\'ebriques des groupes semi--simples complexes.}
Colloque d'alg\`ebre sup\'erieure, Centre Belge de Recherches Math\'ematiques,
\'Etabli\-ssement Ceuterick, Louvain, 1957, 261--289.

\smallskip

[TV] B.~To\"en, M.~Vaqui\'e. {\it Au--dessous de $Spec\,\bold{Z}$.} e-print math.AG/0509684

\smallskip

[van D] D.~van Danzig. {\it Nombres universels ou $\nu !$--adiques avec une introduction
sur l'alg\`ebre topologique.} Ann. Sci. de l' \'ENS, $3^e$ s\'erie, tome 53 (1936), 275--307.

\smallskip

[Za]  D.~Zagier. {\it Vassiliev invariants and a strange identity related to the Dedekind eta--function.}
Topology, 40 (2001), 945--960.

\enddocument